\documentclass[a4paper,14pt]{article}
\usepackage{pstricks,pstcol,pst-grad,pst-node}
\usepackage{latexsym}
\usepackage{amsmath}
\usepackage{amsfonts}
\usepackage{pstcol}
\usepackage{amscd}
\usepackage[all]{xy}
\newcommand{\B}[1]{\mathbb #1}
\newcommand{\C}[1]{{\cal #1}}

\newcommand{\p}{{\partial}}
\newcommand{\al}{{\alpha}}

\newcommand{\be}{{\beta}}

\newcommand{\Om}{{\Omega}}
\newcommand{\om}{{\omega}}
\newcommand{\eps}{{\varepsilon}}

\newcommand{\De}{{\Delta}}
\newcommand{\ga}{{\gamma}}
\newcommand{\Ga}{{\Gamma}}

\newcommand{\Si}{{\Sigma}}
\newcommand{\del}{\partial}

\newcommand{\qed}{\rightline {$\Box $}}

\newcommand{\map}[1]{\stackrel {#1}\longrightarrow}

\newcommand{\Mo}{(M,\omega )}

\newcommand{\pis}{\pi _1(Symp\Mo )}

\newcommand{\im}{\text {im \!}}
\newcommand{\rank}{{\hbox {\em rank }}}

\newcommand{\blowMV}{\widetilde {M}_V}
\newcommand{\cp}{{\B C\B P}}
\newcommand{\sfib}{\Mo \map{i} P\map{\pi } B}
\newcommand{\T }{\hbox {\text {T}}}
\newcommand{\PD }{\hbox {\text {PD}}}

\newcommand{\pf}{\NI {\bf Proof: }}
\newcommand{\mtimes}{\times \dots \times}
\newcommand{\moplus}{\oplus \dots \oplus}
\newcommand{\QED}{\hfill$\Box$\medskip}

\newcommand{\BS}{{\bigskip}}
\newcommand{\NI}{{\noindent}}

\newtheorem{theorem}{Theorem}[section]
\newtheorem{thm}[theorem]{Theorem}
\newtheorem{cor}[theorem]{Corollary}
\newtheorem{defin}[theorem]{Definition}
\newtheorem{rem}[theorem]{Remark}
\newtheorem{lemma}[theorem]{Lemma}
\newtheorem{prop}[theorem]{Proposition}
\newtheorem{ex}[theorem]{Example}
\numberwithin{equation}{section}

\title{Restrictions on symplectic fibrations}
\author{Jaros{\l}aw K\c edra
               \thanks{Partially supported by 
State Committee for Scientific Research grant 2PO3A 035 14. The paper was
written during the author's stay at MPI in Bonn
and at IMPAN in Warsaw.
        \newline
        AMS classification(2000): Primary 57R17;
        Secondary 53D45}\\
 University of Szczecin\\
(with an appendix written jointly with Kaoru Ono)}

\date{\today}
\begin{document}

%USTAWIENIA DO RYSUNKOW
\psset{unit=1.5cm}

\maketitle

\BS

\section{Introduction}\label{S:intro}

This paper is devoted to restrictions on symplectic fibrations
coming from Gromov-Witten invariants. They might have two types
of nature. First is reflected in the properties of the Leray-Serre 
spectral sequence. The behavior of its differentials 
gives an information about the topology of symplectomorphism
groups. There arises two conjectures:

\begin{enumerate}
\item {\bf The flux conjecture} which says that the group
$Ham\Mo $ of Hamiltonian symplectomorphisms is $C^1$-closed
in the group $Symp\Mo $ of all symplectomorphisms, where
$\Mo $ is a compact symplectic manifold. This conjecture
can be expressed equivalently in terms of so called
{\bf flux groups}, which are defined with the use of
the differentials in spectral sequences associated to
symplectic fibrations over two dimensional sphere \cite{lmp1,ms1}. 

\item {\bf The c-splitting conjecture} which states that
the spectral sequence associated to any Hamiltonian fibration
degenerates at $E_2$ term \cite{lm}. This means that the
rational cohomology of the total space is additively isomorphic 
to the tensor product of the cohomology of the base and the
cohomology of the fiber.
\end{enumerate}

The second nature of restrictions arises under assumption
that the spectral sequence collapses at $E_2$, e.g. for
symplectic manifolds satisfying the c-splitting conjecture.
Namely, one may ask how much the cohomology ring of the total
space differs from the tensor product of the cohomology ring
of the base and of the fiber. This was already taken up by Seidel
\cite{se} in case when $\Mo $ is a product of complex
projective spaces.

\subsection{Sketch of the argument}

We are mainly interested in the properties of
the Leray-Serre spectral sequences associated to symplectic
and Hamiltonian fibrations. The main idea comes from the
fact that Gromov-Witten invariants may be defined parametrically.
Here we give the rough sketch of the main points. We need to
slightly restrict the form of fibrations we are interested in.

\begin{defin}\label{D:simple}
Let $M \to P \to B$ be a  fibration whose local coefficient
system is trivial and the base is an $m$-dimensional finite CW-complex
with one dimensional top (co-) homology. 
We call such a fibration {\bf simple}. For example,
$B$ might be a simply connected compact manifold.
\end{defin}

Let $\sfib $ be a symplectic fibration which is simple.
(We will indicate where the simplicity of the fibration is
needed.)
Then 
each fiber $M_b$, $b\in B$, admits a tamed almost 
complex structure $J_b$. These
structures form a smooth family parameterized by the base
of the fibration (see Proposition \ref{P:acs}).
Over each fiber, one has a moduli space 
$\C M_{g,k}(A,\B J_b)$
of 
$\B J_b$-holomorphic stable maps representing a homology
class $A\in H_2(M;\B Z)$ such that $i_*(A)\neq 0$.
Here we use the assumption that the local 
coefficient system of the fibration is trivial
in order to make sense that a map to a fiber $M_b$
represents the class $A\in H_2(M;\B Z)$.
These moduli spaces match up to a bigger moduli space denoted
by $\C M^P_{g,k}(A,\B J)$.

Let $ev_k^P:\C M_{g,k}^P(A, \B J)\to P\oplus \dots \oplus P$
be the evaluation of stable maps at $k$ marked points.
Out of the image of
this mapping we obtain a cycle whose homology class
$[ev_k(\C M_{g,k}(A, \B J))]$
serves for defining a parametrized Gromov-Witten invariant.
Observe that the intersection of this image with
the fiber $M_b$ is the image of the evaluation map
$ev_k:\C M_{g,k}(A, \B J_b)\to M_b\times \dots \times M_b$
from the moduli space associated to this fiber
(see Section \ref{S:gw}).

Here we touch the main point (Theorem \ref{T:main}).
Let us assume that $E_2=E_3=\dots =E_m$ in the spectral
sequence associated to the fibration. It means that
all the differentials are trivial except possibly those
on $E_m$. 
It follows from
the above discussion that the homology class 
$[ev_k(\C M_{g,k}(A,\B J_b))]$
made of the image of the evaluation map for
a generic fiber is the homological intersection
of the homology class coming from the parametrized
moduli space with the fiber of the appropriate
Whitney sum. This implies that the image of
$[ev_k(\C M_{g,k}(A,\B J_b))]$ 
under the differential in the Leray-Serre
spectral sequence is zero. 

A cohomological
version of this argument gives a formula which connects 
Gromov-Witten invariant $\Phi _A $ and the differential 
$\partial $ of the spectral sequence 
(Theorem \ref{T:maincoho}):

$$\sum _{i=1}^k (-1)^{\sum _{1\leq j < i}\deg (\al_j)} 
\Phi _A (\al _1, \dots ,\partial \al_i, \dots ,\al_k)=0 $$
 (see Section \ref{S:main}
for more precise statements and details).

\subsection{A consequence for flux groups}

We present two main applications of the above argument. The first
consists of an estimate of the rank of flux groups. This is
a completion of Theorem A from \cite{k}. We need to introduce
some notions before giving the statement of the result.
Let $ev:Symp\Mo \to M$ denotes the evaluation map,
%%%%%%%%%%%%%%%%%%%%%%%%%%%%%%%%%%%%%%%%%%%%%%%%%%%%%%%%%%
\footnote{ {\bf Warning:} Notice that we use two completely
different evaluation maps. }
%%%%%%%%%%%%%%%%%%%%%%%%%%%%%%%%%%%%%%%%%%%%%%%%%%%%%%%%%
 that is
$ev(\phi )=\phi (pt)$, where $pt\in M$ is a base point.

\begin{defin}\label{D:got}
The image of
$ev_*:\pi_1(Symp\Mo )\to \pi_1(M)$
is called
{\bf symplectic Gottlieb group} and denote
by $G \Mo $. Similarly, the image of
$ev_*:H_1(Symp\Mo ;\B Q)\to H_1(M;\B Q)$
we call a
{\bf rational symplectic Gottlieb group} and denote
by $G_{Q}\Mo $. 
\end{defin}

Recall that the definition of the classical Gottlieb
group uses the space of homotopy equivalences instead of
the group of symplectomorphisms \cite{got}. 
Of course symplectic Gottlieb group is a subgroup of
the classical one.
It is also
worth noticing that the nontriviality of 
$G_Q(M,\om )$ is quite restrictive. Indeed, if this holds
then the rational cohomology of $M$ splits as 
$H^*(M\slash S^1)\otimes H^*(S^1)$ \cite{lm}. Moreover,
elements of Gottlieb group acts trivially on homotopy
groups.

Our method relies on the examination of the $J$-holomorphic
curves in  total spaces symplectic fibrations. 
That's why the homology classes
represented by these curves have to belong to the image
of the map induced by the inclusion of fibers.
This is the motivation for the following notion.

\begin{defin}\label{D:fluxfree}
A homology class $A\in H_2(M;\B Z)$ is said to be
{\bf flux free} if for every $\xi \in \pis $,
$i_*A\neq 0$. Here $\Mo \map{i}P_{\xi }\to S^2$ is the 
symplectic fibration
associated to $\xi $ by the clutching construction.
\end{defin}

%%%%%%%%%%%%%%%%%%%%%%%%%%%%%%%%%%%%%%%%%%%%%%%%%%%%%%%%%%%%%%%%%%%%%%%%%%%%%%
\begin{thm}\label{T:rank}
Let $\Mo $ be a compact symplectic manifold. 
Then the
rank of its flux group $\Ga _{\om }$ satisfies the following
estimate:

$$\rank G\Mo \leq \rank \Ga _{\om }\leq \dim G_{Q}\Mo + 
\dim \big 
[\ker (\cup [\om ]^{n-1}) \cap \ker (\cup  PDA)
\big ],$$
for any flux free class $A\in H_2(M;\B Z)$ such that
$\dim \C M_{g,1}(A,J)=2$. 
\end{thm}
%%%%%%%%%%%%%%%%%%%%%%%%%%%%%%%%%%%%%%%%%%%%%%%%%%%%%%%%%%%%%%%%%%%%%%%%%%%%%%%%

Usually, the above estimate
reduces to a simpler one, which does not involve the
Gottlieb group because it is trivial. 
Moreover, in this case (cf. Proposition \ref{P:flux})

$$
\Ga _{\om }\subset  
\ker (\cup [\om ]^{n-1}) \cap \ker (\cup PDA )
.$$

\NI
Recall that the flux conjecture depends only on the behavior
of the flux homomorphism on the loops with trivial 
(in $\pi_1(M)$) evaluation \cite{lmp1} (Proposition 1.2).
Namely, if the flux group is not discrete then 
its nondicreteness occur in the above intersection of
the kernels. This argument easy shows that the conjecture
holds for manifolds with trivial $\ker \cup [\om ]^{n-1}$,
e.g. K\"ahler manifolds or more generally Lefschetz.

\BS

\subsection{A consequence for Hamiltonian fibrations}

The second application establishes the c-splitting conjecture
for manifolds whose Gromov-Witten invariants satisfy some
restrictions.

%%%%%%%%%%%%%%%%%%%%%%%%%%%%%%%%%%%%%%%%%%%%%%%%%%%%%%%%%%%%%%%%%%%%%%%
\begin{thm}\label{T:ham}
Let $\sfib $ be a compact Hamiltonian fibration. 
Suppose that for any
$\al \in H^i(M)$, where $i \leq 2n -4$,
there exist $\be _0,\be _1,\dots ,\be_k \in \im i^*$  such
that 

$$\Phi _A(\al \cup \be_0,\be_1,\dots,\be_k)\neq 0$$

\NI
for some Gromov-Witten invariant $\Phi _A $.
Then the spectral sequence associated 
to $\sfib $ collapses at $E_2$.
\end{thm}
%%%%%%%%%%%%%%%%%%%%%%%%%%%%%%%%%%%%%%%%%%%%%%%%%%%%%%%%%%%%%%%%%%%%%%%%

\begin{rem}\label{R:ham}{\em \hfill

\begin{enumerate}
\item
One can choose  the elements $\be_i$ from the subalgebra
of $H^*(M)$ generated by the class of symplectic form,
Chern classes and elements of degree at most 2 (or 3 provided
that $B$ is a compact manifold). Indeed, they
are contained in $\im i^*$ as follows. The class
of symplectic form is in this image because the fibration is
Hamiltonian. The Chern classes of $\Mo $ are images of the Chern
classes of the vector bundle tangent to the fibers \cite{k}.
An argument for the elements of small degree is explained in
the proof of the theorem (Section \ref{SS:hf}).

\item
If $B$ is a manifold then the Poincare duality for $P$
allows us to consider elements $\al $ such that
$1\leq \deg (\al )< \frac{2n-b}{2} + 1$, where $2n=\dim M$
and $b=\dim B$. The second inequality is an immediate consequence
of Poincar\' e duality for $P$ and the first one follows from the
fact that $[\om ]^n$ survives since the fibration is Hamiltonian.
\item
By the above observations,
the theorem seems to be relevant in dimension four.
But then c-splitting conjecture easy follows from the theorem
which establishes the conjecture for Hamiltonian fibrations
over 3-dimensional CW-complexes \cite{lm} (Lemma 4.1.8). 
The same argument
works also for simply connected 6-dimensional $\Mo $.
\end{enumerate}
}
\end{rem}

A more explicit application  is the
establishing c-splitting conjecture for projective space
blown-up along 4-dimensional submanifold (Section 
\ref{SS:ex}).

\BS

\noindent
{\bf Acknowledgments.} I am deeply grateful to Kaoru Ono
who patiently answered my questions and pointed out some
gaps in my arguments. I thank Dusa McDuff and Tomek Maszczyk
for valuable comments and remarks.
Also I thank the Max-Planck-Institute,
for providing a wonderful research atmosphere. 

\BS

\section{Symplectic fibrations}\label{S:sf}

\begin{defin}\label{D:sf}
Let $\Mo $ be a symplectic manifold. A fibration
$\Mo \to P \to B$ is said to be {\bf symplectic} if
its structure group is contained in the group $Symp\Mo $
of symplectomorphisms of $\Mo $.
\end{defin}

Recall that an almost complex structure $J$ on $M$ is called
$\om $-tamed if $\om (J\cdot ,\cdot )$ defines
a Riemannian metric. Every symplectic manifold admits
a contractible nonempty set of $\om $-tamed almost complex
structures. This fact can be extended to fibrations
as follows.

\begin{prop}\label{P:acs}
Let $\sfib $ be a symplectic fibration. 
Then there exists a complex bundle $(Vert,\B J) \to P$ such
that $i_b^*\B J$ is $\om _b$-tamed for every $b\in B$.
The set of such almost complex structures $\B J$ is
contractible.
\end{prop}

\pf
We construct a bundle $TM\to Vert\to B$ using the differentials
of the transition functions of the fibration $P\to B$. The
the natural projection $Vert\to P$ defines a 
symplectic vector bundle. This means that 
the structure group of the bundle is $Sp(2n,\B R)$.
Reducing the structure group to the maximal compact
subgroup $U(n)$ we obtain required almost complex
structure. Since such a structure is a section of
the bundle whose fibers are contractible 
($Sp(2n,\B R)\slash U(n)$), then the set of
of such structures is also contractible.

\qed

\section{Parametrized Gromov-Witten invariants}\label{S:gw}

Let $J$ be an $\om $-tamed almost complex structure on $M$ and
$\Si $ be a closed Riemann surface with complex structure $j$.
Let $\C M_{g,k}(A,J)$ be the moduli space of 
stable $J$-holomorphic maps of genus $g$ with $k$ marked
points representing the class $A$. This space is
a compactification of the moduli space of $J$-holomorphic maps 
$u:\Si \to M$ with $k$ marked points, where $\Si $ is a
Riemann surface of genus $g$. The fundamental result of
Fukaya and Ono \cite{fo}(Theorem 1.3)  states that the above
moduli space carries a kind of fundamental class of
$\C M_{g,k}(A,J)$ over $\B Q$ in the following sense.

The $k$-point evaluation map 

$$ev_k:\C M_{g,k}(A,J)\to M\times ... \times M$$ 
is defined by
$$ev_k(u,z_1,...,z_k)= (u(z_1),...,u(z_k)).$$

\NI
Then the result of Fukaya and Ono says that
$[ev_k(\C M_{g,k}(A,J))]$ is well defined as a homology class
in $H^*(M\times \dots \times M;\B Q)$. This allows them to
define a Gromov-Witten
invariant as a map 
$\Phi _A :H^*(M)\times \dots \times H^*(M) \to \B Q$ by

$$\Phi _A (\al_1,...,\al_k) = 
\left <a_1\times ...\times a_k, [ev_k(\C M_{g,k}(A,J))]\right >.$$

The concept of the Gromov-Witten invariants
can be generalized to the parametrized situation.
Namely, let $\sfib $ be a symplectic fibration over
a finite simply connected CW-complex. 
Then every fiber is, up to 
symplectic isotopy, identified with  $\Mo $. 
Suppose that $i_*(A)\neq 0$
for some $A\in H_2(M;\B Z)$. Then
for an almost complex structure $\B J$ in
the vertical bundle of a symplectic fibration,
we consider the space $\C M^P_{g,k}(A,\B J)$ 
of pairs $(u,b)$ where
$b\in B$ and $u:C_g\to P_b$ 
is a genus $g$, $\B J_b$-holomorphic
map representing the class $A$. It does make sense to say that
some map $\Si \to P_b$ represents class $A$ since every
fiber is identified with $M$ up to symplectic isotopy.

A parametrized $k$-point evaluation map is  defined
as follows.

$$
\CD
ev^P_k:\C M^P_{g,k}(A,\B J)\to P\oplus \dots \oplus P\\
ev^P_k(u,b,z_1,...,z_k)=(u(z_1),...,u(z_k)),
\endCD
$$
where $P\oplus \dots \oplus P$ is $k$-fold Whitney sum of the
fibration $P$:

$$
\CD
P\oplus \dots \oplus P:=\De ^*(P\times \dots \times P)
@>\widetilde {\De }>> P\times \dots \times P\\
\qquad\qquad\qquad\qquad @VVV @VVV\\
\qquad\qquad\qquad\qquad B @>\De >> B\times \dots \times B
\endCD
$$

\BS
\NI
To define the parametrized Gromov-Witten invariant
we need to construct a homology class \\
$[ev^P_k(\C M_{g,k}^P(A,\B J))]\in H_*(P\moplus P)$. To do this one has
to observe that the construction of Fukaya and Ono
of the cycle $ev_k(\C M_{g,k}(A,J))$ can be done parametrically.
Finally the parametrized Gromov-Witten invariant
is a map
$\Phi^P_A :H^*(P)\times \dots \times H^*(P) \to \B Q$ defined by

$$\Phi^P _A (\al_1,...,\al_k) = 
\left <\widetilde {\De }^*(\al_1\times \dots\times \al_k) , 
\left [ev_k^P(\C M^P_{g,k}(A,\B J))\right ]\right >.$$

It is obvious from the definition that
$ev_k^P(\C M_{g,k}^P(A,\B J))\cap (P\oplus \dots \oplus P)_b 
=ev_k(\C M_{g,k}(A,\B J_b))$ for any symplectic fibration over
simply connected base. We need an analogous statement in
homology. In more detail, let 

$$
i_!:H_{k+m}(P\oplus \dots \oplus P)\to H_k(M\times \dots \times M)
$$
be the homology transfer, where 
$i:M\times \dots \times M \to P\oplus \dots \oplus P$ is an inclusion
of the fiber and $\dim B=m$. Recall that the transfer can be 
expressed as follows.

$$
i_!(a)= [C_a \cap i(M\times \dots \times M)],
$$
where $C_a$ is a cycle representing the class $a$.

\begin{lemma}\label{L:main}
Let $\B J$ be an almost complex structure in $Vert$
compatible with the fibration.
Then  

$$
i_! \left [ev^P_k(\C M^P_{g,k}(A,\B J)\right ] = 
[ev_k(\C M_{g,k}(A,\B J_b)],$$
for any $b\in B$.

\end{lemma}

\pf
This follows from the above observation that the construction of
the cycle  $ev_k (\C M_{g,k}(A,J)$ in \cite{fo} can
be done parametrically. 
We proceed using the framework of \cite{fo}.

Fix a point $b\in B$ and construct Kuranishi structure
on the moduli space 
$\C M_{g,k}(A,\B J_b)$. 
Next choose
multivalued perturbation in order to get a virtual
fundamental cycle. Then extend the Kuranishi structure to
a Kuranishi structure on the parametrized moduli space
and also extend the multivalued pertubation (for
$\C M_{g,k}(A,\B J_b)$)  to a multivalued pertubation
of the Kuranishi map for the parametrized moduli space.

For the above choice of the Kuranishi structure and
the multivalued perturbation, the intersection 
$(P\moplus P)_b\cap ev^P_k(\C M^P_{g,k}(A,\B J))= 
ev_k(\C M_{g,k}(A,\B J_b)$ is transversal, which
completes the proof.
\QED

\begin{rem}\label{R:main}
{\em
A similar argument should also work for the other approaches
to Gromov-Witten invariants, such as those presented by
Li and Tian \cite{lt} or Siebert \cite{si}. 
The case of curves of genus zero was
already done by Seidel \cite{se}. Notice also that in the
weakly monotone case the above lemma might be proved using
genericity of almost complex structure involved.}
\end{rem}

\begin{ex}\label{E:hir}
{\em
Let $\Mo:=(S^2\times S^2,\om )$ be equipped with symplectic structure
such that $[\{x\}\times S^2]$ has area 1 and $[S^2\times \{x\}]$
has area 2. Let $A:=[S^2\times \{x\}] + [\{x\}\times S^2]$
There exist a symplectic fibration

$$\Mo \to P\to S^2$$

\NI
which admits a compatible almost complex structure $\B J$ in $Vert$
such that $\C M_{0,3}(A,\B J_z)$ is empty for any $z\in S^2$ except
one point. 
This phenomenon can be explained as follows (cf. \cite{am}). 
The space of 
$\om $-tamed almost complex structures contains an open
stratum consisting almost complex structures $J$ for which
class $A$ cannot be represented by any $J$-holomorphic
sphere and codimension 2 stratum consisting of almost
complex structures for which $A$ can be represented.
The complex structure $\B J$ can be seen as a map
from $S^2$ to the space of $\om $-tamed almost complex
structures which intersects transversely this codimension
2 stratum at one point.\QED
}
\end{ex}

\section{The main theorem}\label{S:main}

Our main results gives the restrictions to the form of the
differentials in the 
Leray-Serre spectral sequences associated to symplectic
fibrations. We consider the spectral sequences for
fibration $P\to B$ as well as for $P\moplus P\to B$.
Corresponding differentials in the homology spectral
sequences we denote by $\partial ^r:E^r_{p,q}\to E^r_{p-r,q+r-1}$
and $\overline \partial ^r:\overline {E}^r_{p,q}\to 
\overline E^r_{p-r,q+r-1}$, respectively. The analogous
differential in the cohomology spectral sequence we
will denote with subscripts.

\subsection{Generalized Wang homomorphisms}\label{SS:wang}

Now we define a generalizations of  the Wang homomorphisms,
which originally were associated to fibrations over homology
spheres (cf. \cite{sp} Section 9.5). Suppose
that the differentials in the Leray-Serre cohomology spectral
sequence associated to a simple fibration over 
$m$-dimensional base
are trivial up to the $m^{\text {th}}$ term.
In other words $\partial ^m$  (respectively $\partial _m$) 
is the only nontrivial differential in the homology
(resp. cohomology) spectral sequence. 
For the differential 
$\partial _m :H^0(B)\otimes H^k(M) \to H^m(B)\otimes H^{k-m+1}(M)$
we define a {\bf generalized Wang homomorphism} 
$\partial :H^k(M)\to H^{k-m+1}(M)$ by 

$$\partial _m(1\otimes \al ) = \be \otimes \partial \al,$$
where $\be \in H^m(B)$ is a generator dual to the fundamental
class $[B]$. 

\begin{lemma}\label{L:wang}

The generalized Wang homomorphism:
\begin{enumerate}
\item
is natural with respect to bundle maps which induce
the identity on the top cohomology of the base and
\item
satisfies the Leibniz rule.
\end{enumerate}
\end{lemma}

\pf We use the fact that the differential in the spectral
sequence has required properties.

\begin{enumerate}

\item Naturality:

\noindent
Let $f:P_1\to P_2$ be a bundle map.

\begin{eqnarray*}
&&\be \otimes \partial (f|M)^*\al =\\
&&\partial _m(1 \otimes (f|_M)^*\al )=\\
&&\partial _m(f^*(1 \otimes \al))=\\
&&f^*(\partial _m(1 \otimes \al))=\\
&&f^*(\be \otimes \partial \al)=\\
&&\be \otimes (f|M)^*\partial \al)
\end{eqnarray*}

\item The Leibniz rule:
\begin{eqnarray*}
&&\be \otimes \partial (\al _1 \cup \al _2)=\\
&&\partial _m(1\otimes (\al _1 \cup \al _2)=\\
&&\partial _m((1\otimes \al _1)\cup (1\otimes \al _2))=\\
&&\partial _m(1\otimes \al _1) \cup (1\otimes \al _2) \pm
(1\otimes \al _1)\cup \partial _m(1\otimes \al _2)=\\
&&\be \otimes \partial (\al _1) \cup (1\otimes \al _2) \pm
\cup (1\otimes \al _1) \cup \be \otimes \partial (\al _2)=\\
&&\be \otimes (\partial (\al _1) \cup \al _2 \pm
\al _1 \cup \partial (\al _2)).
\end{eqnarray*}

\end{enumerate}
\qed

The homology case, denoted by the same symbol,
$\partial :H_{k-m+1}(M)\to H_k(M)$ (that should not be confusing)
is defined analogously by

$$\partial ^m([B]\otimes a) = [pt]\otimes \partial a.$$
It is also natural with respect to bundle maps which
induce the identity on the top homology of the base.
Moreover, it straightforward that homology and cohomology Wang homomorphisms
are dual to each other in the sense that

$$\left <\partial a,\al\right > = 
\left <a , \partial \al \right >.$$

\subsection{The Wang homomorphism vanishes on Gromov-Witten invariants}
\label{SS:main}

Let $\overline \partial $ denote
the generalized Wang homomorphisms associated to the Whitney
sum of the fibration.

%%%%%%%%%%%%%%%%%%%%%%%%%%%%%%%%%%%%%%%%%%%%%%%%%%%%%%%%%%%

\begin{thm}\label{T:main}
Let $\sfib $ be a simple (in the sense of Definition
\ref{D:simple})
symplectic fibration over $m$-dimensional 
base and
$A\in H_2(M;\B Z)$ be such that $i_*A\neq 0$.
Suppose that the differential
$\partial ^r :E^r_{*,*}\to E^r_{*-r,*+r-1}$  
is zero for $r\neq m $.
Then 
$\overline \partial [ev_k(\C M_{g,k}(A,J))]  = 0$
for any moduli space $\C M_{g,k}(A,J)$.
\end{thm}

%%%%%%%%%%%%%%%%%%%%%%%%%%%%%%%%%%%%%%%%%%%%%%%%%%%%%%%%%%%

\pf
Since $E^m_{*,*}$ is the only term with possibly 
nontrivial differential then so does $\overline E^m_{*,*}$.
Thus
there is the following commutative diagram in which the
upper row is exact 

$$
\!\!\!\!\!\!\!
\CD
{\scriptstyle
\!\!0\to \overline E^{\infty }_{m,k}} @>>> 
{\scriptstyle \overline E^m_{m,k}} @> {\overline \partial  ^m}>>
{\scriptstyle \overline E^m_{0,k+m-1} }
&\to &
{\scriptstyle \overline E^{\infty }_{0,k+m-1}\to 0}\\
@AAA @| @| @|\\
{\scriptstyle
H_{m+k}(P\moplus P)} @>{i_!}>>
{\scriptstyle  H_k(M \mtimes M)} @>\overline \partial >>
{\scriptstyle  H_{k+m-1}(M\mtimes M)}
&\to &
{\scriptstyle i_*H_{k+m-1}(M \mtimes M)\to 0}
\endCD
$$

\BS
\NI
Here $H_{m+k}(P\oplus \dots \oplus P)\to 
\overline E^{\infty }_{m,k} =
H_{m+k}(P\oplus \dots \oplus P)\slash 
F_{m-1}(H_{m+k}(P\oplus \dots \oplus P)$ is the projection
and $F_s$ denotes the filtration to which is associated 
$\overline E^{\infty }$ (see \cite{sp}).
It follows from Lemma \ref{L:main}, that 

$$
i_! [ev^P_k (\C M^P_{g,k}(A,\B J))] =
[ev_k(\C M_{g,k}(A,J))].
$$

\NI
Then, accordingly to the exactness of the above sequence,
we get that 

$$
\overline \partial [ev_k(\C M_{g,k}(A,J))]=0
$$

\NI
as required. \QED

\begin{rem}
{\em
Of course, we can relax the assumption that there is the only
one possibly nontrivial differential. Indeed,
since the class $[\C M_{g,k}(A,J)]= i_![\C M^P_{g,k}(A,\B J)]$
then it is in the kernel of any differential in the spectral
sequence.
The motivation for this assumption is the following. In case of
Hamiltonian fibrations, it is already proven that the first
nontrivial term in the spectral sequence is $E_4$. By taking
appropriate restrictions we can prove c-splitting in some
cases by induction.
On the other hand, in non Hamiltonian case we are mostly
interested in fibrations over $S^2$.}
\end{rem}

\subsection{The Wang homomorphism is compatible with Gromov-Witten
invariants}

The following theorem is a cohomological incarnation of
the previous one.

%%%%%%%%%%%%%%%%%%%%%%%%%%%%%%%%%%%%%%%%%%%%%%%%%%%%%%%%%%%%%%%%%%%%%%%%%%%%

\begin{thm}\label{T:maincoho}
Under the assumption of Theorem \ref{T:main}, the following
symmetry holds

$$\sum _{i=1}^k (-1)^{\sum _{1\leq j < i}\deg (\al_j)} 
\Phi _A 
(\al _1, \dots ,\partial \al_i, \dots ,\al_k)=0, $$
where $\al_i\in H^*(M)$. In particular,

$$\Phi_A (\partial \al, \al_1,\dots ,\al_k)=0$$
for $\al_i\in \ker \partial$.
\end{thm}

%%%%%%%%%%%%%%%%%%%%%%%%%%%%%%%%%%%%%%%%%%%%%%%%%%%%%%%%%%%%%%%%%%%%%%%%%%%%

\pf This is the following computation, which uses the 
naturality of the generalized Wang homomorphism and the fact that
it satisfies the Leibniz rule.

\begin{eqnarray*}
&&\sum _{i=1}^k (-1)^{\sum _{1\leq j < i}\deg (\al_j)} 
\Phi_A 
(\al _1, \dots ,\partial \al_i, \dots , \al_k)=\\
&&
\sum _{i=1}^k (-1)^{\sum _{1\leq j < i}\deg (\al_j)} 
\left <
\al _1 \times \dots \times \partial \al_i\times \dots \times \al_k,
[ev(\C M_{g,k}(A,J))]\right >=\\
&&
\sum _{i=1}^k (-1)^{\sum _{1\leq j < i}\deg (\al_j)} 
\left <
\pi _1^*\al _1 \cup \dots \cup \pi _i^* \partial \al_i\cup \dots \cup \pi _k^* \al_k,
[ev(\C M_{g,k}(A,J))]\right >=\\
&&
\sum _{i=1}^k (-1)^{\sum _{1\leq j < i}\deg (\al_j)} 
\left <
\pi _1^*\al _1 \cup \dots \cup 
\overline {\partial }\pi _i^* \al_i\cup \dots \cup \pi _k^* \al_k,
[ev(\C M_{g,k}(A,J))]\right >=\\
&&\left <\overline{\partial } (\al_1\times \dots \times \al_k),
[ev(\C M_{g,k}(A,J))]\right >=\\
%&&\\
%
%
&&\left <\al_1\times \dots \times \al_k,
\overline {\partial } [ev(\C M_{g,k}(A,J))]\right >=0
\end{eqnarray*}
\QED

\section{Consequences}\label{S:con}

It follows from the work of Lalonde, McDuff and Polterovich 
\cite{lmp1,m2} that the spectral sequence for Hamiltonian
fibration over $S^2$ degenerates at $E_2$ term. Moreover,
the same holds for Hamiltonian fibrations over any 3-dimensional
CW-complex  and it is conjectured by Lalonde and McDuff
to be  true in general \cite{lm}.
We divide the consequences of the main theorem into two parts.
First, symplectic fibrations over $S^2$ which are not
Hamiltonian. Second, Hamiltonian fibrations over arbitrary
bases.

\subsection{Fibrations over $S^2$ and flux groups}\label{SS:s2}

Every symplectic fibration over $S^2$ is naturally
associated to an element $\xi \in \pi _1(Symp\Mo )$
by the clutching construction. We denote this fibration
by $P_{\xi }$. Recall that the evaluation map
$ev:Symp\Mo \to M$ is defined by $ev(\phi )=\phi (pt)$,
where $pt\in M$ is a point.

%%%%%%%%%%%%%%%%%%%%%%%%%%%%%%%%%%%%%%%%%%%%%%%%%%%%%%%%%%%%%%%%%%%%
\begin{prop}\label{P:evn0}
Let $\Mo \to P_{\xi } \to S^2$ be a symplectic fibration
and $A \in H_2(M;\B Z)$ is such that $i_*A\neq 0$.
Suppose that $[ev_*(\xi )]\neq 0$ in homology.
Then 
$$\Phi _A (\partial \al_1,\dots ,\partial \al_k)=0,$$
for any Gromov-Witten invariant.
\end{prop}
%%%%%%%%%%%%%%%%%%%%%%%%%%%%%%%%%%%%%%%%%%%%%%%%%%%%%%%%%%%%%%%%%%%%

\pf 
It is a result of Lalonde and McDuff \cite{lm}, that
$[ev_*(\xi )]\neq 0$ in homology iff 
$\ker \partial = \im \partial$. Thus the statement easily
follows from Theorem \ref{T:maincoho}.
\QED

The next result gives the restriction for the rank of flux
groups. We recall the definition of these groups.
{\bf Flux homomorphism} $F:\pi _1(Symp\Mo )\to H^1(M;\B R)$
is defined by

$$F(\xi ) = \partial _{\xi }[\om ],$$
where $\partial _{\xi }$ is the Wang homomorphism
associated to the fibration $P_{\xi }$. By definition
{\bf the flux group } $\Ga _{\om }$ is the image
of the flux homomorphism. The importance of flux groups
comes from the fact that they, in some sense, measure
the difference between the group of symplectomorphisms
$Symp\Mo $
and the group of Hamiltonian symplectomorphisms
$Ham \Mo $. Moreover,
the discreteness of flux groups in equivalent to closeness
of $Ham \Mo $ in $Symp\Mo $ (see \cite{lmp1,ms1} for details).

%%%%%%%%%%%%%%%%%%%%%%%%%%%%%%%%%%%%%%%%%%%%%%%%%%%%%%%%%%%%%%%%%%%%
\begin{prop}\label{P:flux}
Let $\Mo \map{i} P_{\xi } \to S^2$ be a symplectic fibration.
Suppose that $[ev_*(\xi )]=0$ in homology and $F(\xi )\neq 0$.
Then

$$F(\xi )\in \ker (\cup [\om ]^{n-1})
\cap \ker (\cup PDmA)$$
for any  $A\in H_2(M;\B Z)$ for which 
$\dim \C M_{g,1}(A,J)=2$ and $i_*A\neq 0$.
Here $m\in \B Q$ denotes the virtual number of stable
curves representing the class $A$, that is
$mA = [ev_1(\C M_{g,1}(A,J))]$.
\end{prop}
%%%%%%%%%%%%%%%%%%%%%%%%%%%%%%%%%%%%%%%%%%%%%%%%%%%%%%%%%%%%%%%%%%%%%

\pf The statement that $F(\xi)\in \ker (\cup [\om ]^{n-1})$
was already proved in \cite{k}. The proof relies on the fact that
the assumption that $[ev_*(\xi )]=0$ in homology implies
that $\partial _{\xi }[\om ]^n = 0$. The latter is equal
to $n \partial _{\xi }[\om ]\cup [\om ]^{n-1} =
F(\xi )\cup [\om ]^{n-1}$ and the statement follows.

The inclusion 
$F(\xi )\in \ker (\cup  PDmA)$
follows from Theorem \ref{T:main} and the observation
that $mA = [ev_1(\C M_{g,1}(A,J))]$. Suppose that 
$m\neq 0$, otherwise the statement is trivial. 
Let $\al \in H^1(M)$ be an arbitrary element. Then
we compute that

\begin{eqnarray*}
&&\Big <F(\xi )\cup PD mA, PD \al \Big>=\\
&&\Big <F(\xi )\cup \al,mA\Big >=\\
&&\Big <\partial _{\xi } [\om ]\cup \al,[ev_1 (\C M_{g,1}(A,J))]\Big >=\\
&&\Big <\partial _{\xi }([\om ]\cup \al ),[ev_1 (\C M_{g,1}(A,J))]\Big >=\\
&&\Big <[\om ]\cup \al,\partial _{\xi }[ev_1 (\C M_{g,1}(A,J))]\Big >=0.
\end{eqnarray*}

\noindent
The assumption on $[ev_*(\xi )]$ implies that the Wang
homomorphism is trivial on $H^1(M)$ which gives the
third equality. 

\qed

\BS

\begin{rem}
The  other assumptions which imply
the hypothesis of the above theorem are discussed in \cite{k}. 
More informations about the
flux groups can be found in \cite{lmp1}.
\end{rem}

\subsection{A proof of Theorem \ref{T:rank}}\label{SS:rank}

For convenience of the reader we recall the statement:

%%%%%%%%%%%%%%%%%%%%%%%%%%%%%%%%%%%%%%%%%%%%%%%%%%%%%%%%%%%%%%%%%%
\BS
\NI
{\bf Theorem \ref{T:rank}}
{\em
Let $\Mo $ be a compact symplectic manifold. 
Then the
rank of its flux group $\Ga _{\om }$ satisfies the following
estimate:

\medskip
$\rank G\Mo \leq \rank \Ga _{\om }\leq \dim G_{Q}\Mo + 
\dim \big 
[\ker (\cup [\om ]^{n-1}) \cap \ker (\cup PDA)
\big ],$

\medskip
\NI
for any flux free class $A\in H_2(M;\B Z)$
such that
$\dim \C M_{g,1}(A,J)=2$. 
}
%%%%%%%%%%%%%%%%%%%%%%%%%%%%%%%%%%%%%%%%%%%%%%%%%%%%%%%%%%%%%%%%%%%%
\BS

\NI
\pf
We use the commutativity of the following diagram:

$$
\xymatrix@1{
\pi _1(Symp\Mo ) \ar[r]^-{ev_*} \ar[dr]^F & 
                         \pi _1(M) \ar[r] & 
                 H_1(M;\B Z) \ar[r]^-{PD} & 
          H^{2n-1}(M;\B Z)\ar[d]^{\cdot n}\\
                                          &
   H^1(M;\B R) \ar[rr]^{\cup[\om ]^{n-1}} & 
                                          &
                         H^{2n-1}(M;\B R)                                       
%
%\ar[r]-^{ev_*} &  \ar[r] & 
%\ar[d]\\
%& H^1(M;\B R) \ar[r]-^{\cup [\om ]^{n-1}} & & H^{2n-1}(M;\B R)
}
$$
\BS

\noindent
which is easy to prove. Let us consider the first inequality.
Take elements 
$\xi _1, \dots , \xi _k$ ($k=\rank G\Mo $) such that
$\sum a_iF(\xi _i) = 0$ for some nonzero $a_i\in \B Z$.
This means that the element $\sum a_i\xi _i\in \pi _1(Symp\Mo )$
may be represented by a Hamiltonian loop. Hence their
evaluation $ev_*(\sum a_i\xi _i )= \sum a_i ev_*(\xi _i)=0$
due to results from Floer theory and we obtain the first inequality. 
It may be confusing that we use addition in $\pi _1(M)$, 
but it is in fact in the image of $\pi _1(Symp\Mo )$ which is abelian.

The proof of the second inequality requires more advanced
argument which relies on so called {\bf topological rigidity of
Hamiltonian loops} \cite{lmp2}:

\BS
{\em 
If $\xi \in \pi _1(Diff(M))$ can be represented by loops
in $\pi _1(Symp\Mo )$ and $\pi _1(Symp(M,\om ^{\prime }))$,
then $\partial _{\xi}[\om ]=0$ iff $\partial _{\xi}[\om^{\prime} ]=0$}

\BS

\noindent
Let $K:=\ker ev_*:H_1(Symp\Mo ;\B Q)\to H_1(M;\B Q)$. Then
we obtain the extension 
$0\to F(K)\to \Ga _{\om }\to \Ga _{\om }\slash F(K)$,
where $\Ga _{\om }\slash F(K)$ has no torsion. Indeed,
if $F(\xi )\cdot F(K)$ were of finite order, say $k$ then it would mean
that $0=ev_*(k\xi )=kev_*(\xi )$ but $ev_*(\xi )\neq 0$ in 
$H_1(M;\B Q)$ which is impossible. Hence we get that 

$$\text {rank} \Ga _{\om } = 
\text {rank} F(K) + \text {rank} (\Ga _{\om }\slash F(K)).$$

\NI
It follows from the above diagram that
$\text {rank} (\Ga _{\om }\slash F(K))=\text {rank} G_Q\Mo .$ Thus we have only
to show that 

$$\text {rank} F(K) \leq 
\dim \ker (\cup [\om ]^{n-1}) \cap \ker (\cup PDA).$$

\NI
The argument is similar to that in \cite{lmp2} (Theorem 2.D). 
Let $\om _{\eps }$ be a rational symplectic form which is
a small perturbation of $\om $. Then obviously
$\dim \ker \cup [\om]^{n-1}\leq \dim \ker \cup [\om _{\eps }]^{n-1}.$
Moreover, the Gromov-Witten invariant $[ev_1(\C M_{g,1})]$ is
the same for both symplectic forms.

Suppose that there exist elements $\xi _1,...,\xi _l \in \pis$
such that $F(\xi _1),..., F(\xi _l) \in F(K)$
are linearly independent over $\B Z$ and
$l>$ dim $\ker\cup [\om ^{n-1}]\cap \ker \cup PDA$.
Let
$\xi _1^{\epsilon },...,\xi _l^{\epsilon }\in \pi _1(Symp(M,\om _\epsilon )$
are represented by  perturbations of loops representing
$\xi _1,...,\xi _l \in \pis$.
Then the fluxes
$F_{\epsilon }(\xi _j^{\epsilon })$
are rational classes and
are contained in
$\ker\cup[\om _\epsilon ^{n-1}] \cap \ker \cup PDA$, 
according to Proposition \ref{P:flux}.
It follows that some of their
nontrivial combination over $\B Z$ equals zero:
$\sum _j m_jF_{\epsilon }(\xi _j^{\epsilon })=0$
, $m_j\in \B Z$, $j=1,...,k$.
Due to the topological rigidity of Hamiltonian loops, we get that
$\sum _j m_j F(\xi _j)=0$
which gives the contradiction and completes the proof

\qed

\subsection{Hamiltonian fibrations}\label{SS:hf}

The aim of this section is to prove Theorem \ref{T:ham}.
The idea of the proof is to show that
an element $\al \in H^*(M)$ cannot be in the image
of the Wang homomorphism, provided that certain Gromov-Witten
invariant is nonzero. If any $\al \in H^*(M)$ 
has this property, then c-splitting holds for $\Mo $.

\BS
%%%%%%%%%%%%%%%%%%%%%%%%%%%%%%%%%%%%%%%%%%%%%%%%%%%%%%%%%%%%%%%%%%%%%%%

\NI
{\bf Theorem \ref{T:ham}}
{\em
Let $\sfib $ be a compact Hamiltonian fibration. 
Suppose that for any
$\al \in H^i(M)$, where $i \leq 2n -4$,
there exist $\be _0,\be _1,\dots ,\be_k \in \im i^*$ such
that 

$$\Phi _A (\al \cup \be_0,\be_1,\dots,\be_k)\neq 0$$

\NI
for some Gromov-Witten invariant $\Phi _A $.
Then the spectral sequence associated to 
$\sfib $ collapses at $E_2$.
}
%%%%%%%%%%%%%%%%%%%%%%%%%%%%%%%%%%%%%%%%%%%%%%%%%%%%%%%%%%%%%%%%%%%%%%%%

\BS

\pf We use the induction on the dimension of the base. 
Due to Lalonde and McDuff \cite{lm} we know that 
the spectral sequence collapses at $E_2$ for any
Hamiltonian fibration over the base of dimension
at most  3. This implies that $E_2=E_3=E_4$ in the 
spectral sequence associated to an arbitrary fibration.

Suppose the fibration is does not c-split and $E_m$ is the
first nontrivial term in the spectral sequence.
Then there exists some $\eta \in H^*(M)$ such that
$\partial _m(1\otimes \eta) \neq 0$. We restrict
the fibration over an $m$-dimensional CW-complex
with one dimensional top cohomology such that
$\partial _m(1\otimes \eta) \neq 0$ in the induced
spectral sequence.

Let's consider
the associated Wang homomorphism for which we have
that $\partial \eta \neq 0$.
First notice that 
$\deg (\partial \eta ) \leq 2n-4$. 
Indeed, it
follows from the fact that 
the Wang homomorphism coming from the first possible nontrivial
differential decreases the degree by 3 and
$\partial :H^{2n}(M)\to H^{2n-r+1}(M)$
is zero since the top cohomology is generated by $[\om ]^n$.

Since $\be _0,\dots , \be _k\in \ker \del$ then
Theorem \ref{T:maincoho} implies that

$$\Phi _A (\partial \eta \cup \be _0,\be _1,\dots ,\be _k)=
\Phi _A (\partial (\eta \cup \be_0),\be_1,\dots ,\be_k)=0,$$
which contradicts the assumption.

\qed

\BS
As a corollary we can recover the Blanchard theorem which
states that c-splitting holds for manifolds satisfying the
Hard Lefschetz condition.

\begin{cor}[Blanchard \cite{bl}]\label{C:blanchard}
If $\Mo $ satisfies the Hard Lefschetz condition (i.e.
$\cup [\om ]^k :H^{n-k}\to H^{n+k}$ is an isomorphism for
$k=1,...,n$, $\dim M=2n$), then any Hamiltonian fibration
$\sfib $ over a compact manifold $B$ c-splits.
\end{cor}

\pf
Consider the spherical ($g=0$) invariant 
$\Phi _0: H^*(M)\times H^*(M)\times H^*(M)\to \B Q$
for trivial homology class A=0 (see Section \ref{S:gw}). 
Clearly, it is
defined by the usual cup product, namely
$\Phi _0(\al, \be, \ga ) = \left <\al \cup \be \cup \ga ,[M]\right >$
\cite {ms2}.
As usual we restrict ourselves to  fibrations 
with the only one possibly nontrivial term in the spectral
sequence. 
Consider the Wang homomorphism $\del $ in this case. 
Suppose that $\del \eta = \al $, where $\al \in H^k(M)$ and $k< n$
is the least degree of nonzero element in the image of $\del $.
The assumption $k\leq n$ is due to Remark \ref{R:ham} (2).

Let $\be \in H^{2n-k}$ be such that
$\al \cup \be  \neq 0$ in $H^{2n}(M)$.
Then 
$\be = \ga \cup [\om ]^{n-k}$, due to Hard Lefschetz condition, and

$$\Phi _0(\al , [\om ]^k, \ga) =
\be \cup [\om ]^k \cup \ga \neq 0.
$$
Note that $\del (\ga )=0$ because $\deg (\del (\ga )) < k=deg (\al )$.
According to Theorem \ref{T:ham} and 
Remark \ref{R:ham} (2), 
we get the statement.

\qed

\BS

\subsection{An explicit example}\label{SS:ex}

Let $(V, \om _V)$ be a compact 4-dimensional symplectic manifold.
According to results of Gromov \cite{pdr} and
Tischler \cite{ti}, there exists a symplectic
embedding $V\to \cp ^5$. Let's consider a symplectic blow-up
$\Mo :=(\widetilde {\cp ^5}_V,\om )$ of $\cp ^5$ along $V$.
The aim of this section is to prove the following

\begin{prop}\label{P:ex}
Let $B$ be a compact manifold. Then, any Hamiltonian
fibration $\Mo \to P\to B$ is c-split.
\end{prop}

\begin{rem}\label{R:ex}
{\em \hfill

\begin{enumerate}
\item
The c-splitting conjecture is established for symplectic
manifolds satisfying the hard Lefschetz condition (e.g. K\" ahler
manifolds) \cite{bl}. Also, as we already have mentioned,
it is true for 4-dimensional and simply connected 6-dimensional
manifolds. The above proposition proves the c-splitting conjecture
for the family of symplectic manifolds which contains
the easiest examples which don't satisfy  the hard Lefschetz
condition \cite{m1}.
\item
Recall that the cohomology of the above blow-up has the following
form \cite{m1,gi}

$$
H^*(M)=H^*(\cp ^5)\oplus H^*(V)[u]\slash u^3.
$$
\end{enumerate}
}
\end{rem}

\BS

\NI
{\bf Proof of Proposition \ref{P:ex}:}
Throughout the proof we adopt the notation from the
appendix.
Consider the Leray-Serre spectral sequence associated
to the fibration $\Mo \to P\to B$. Due to Lalonde and
McDuff \cite{lm}, we know that the differentials $\del _2 = \del _3 =0$.
Suppose that $\del _4\neq 0$. We can restrict the fibration
over a compact 4-dimensional manifold, such that $\del _4=\del $ is
also nonzero.

\NI
{\bf Claim: $\del :H^5(M)\to H^2(M)$ is zero.}

\NI
Suppose that is there exists $\al \in H^5(M)$
such that $\del (\al ) = ax + bu$, where 
$x=p^*\om _0$, $u:=t^*(\tau )$ and $a,b\in \B R$.
Since the symplectic structure om $M$ away from the
blow-up locus is the same as the standard structure
$\om _0$  on $\cp ^5$,
then

$$
\Phi _L([\om ^5],[\om ^5]) \neq 0,
$$
where $L\in H_2(M)$ is the class of line.
For $a\neq 0$ we get that 
$\Phi _L((ax + bu)\cup {1\over a}x^4, [\om ^5])\neq 0$,
because $u \cup x^4 = 0$. It follows from Theorem \ref{T:maincoho}
that for $a\neq 0$  $ax + bu \notin \im \del $.

Next we show that $u$ cannot be in the image of $\del $.
To see this consider the following Gromov-Witten invariant

$$
\Phi _A(u \cup t^*(\tau _{\om V}),u \cup t^*(\tau _{\om V}))\neq 0.
$$
Since $\del (t^*(\tau _{\om V})=0$ because $M$ is simply connected, 
then we obtain that $u\notin \im \del $. Thus we have proven
the claim.

The rest is easy now. By the Poincar\' e duality, it follows
from the claim that $\del :H^8(M)\to H^5(M)$ is zero.
Notice that $H^8(M)$ is generated by $x^4$ and 
$u \cup t^*(\tau _{\om _V^2})$ hence we get that

$$
\del (t^*(\tau _{\om _V^2})) = 0,
$$
which implies that $\del H^6(M)\to H^3(M)$ is zero.
Again by the Poincar\' e duality we get that
$\del H^7(M)\to H^4(M)$ is zero, which finishes
the proof.

The case of higher dimensional bases go through
in the same way.

\qed

\BS

\appendix
\NI
\section{ Appendix (by Jaros\l aw K\c edra and Kaoru Ono):\\
 Simple examples of nontrivial Gromov-Witten invariants }

\BS

\subsection{The symplectic blow-up construction}

\BS

Let $(M^{2n},\om _M)$ and $(V^{2m},\omega _V)$ be compact
symplectic manifolds such that there exists
a symplectic embedding $(V,\om _V)\to (M,\om _M) $ of codimension
$2k + 2$. We consider a symplectic blow-up of $\Mo $
along $V$ denoted by $\blowMV$.
It is obtained by cutting
out a small tubular neighborhood of $V$ in $M$ and
glue in a symplectic disc bundle over projectivized
normal bundle of $V$. More precisely, there exist a small
tubular neighborhood $\C N_V$ of $V$ in $M$, which is symplectomorphic
to an $\eps $-disc subbundle $\nu _V(\eps )$
in the normal bundle $\nu _V$
of $V$. Since $\nu _V$ is a complex bundle then
then we have that
$\nu _V(\eps )=P\times _{U(k+1)}D^{2k+2}_{\eps }$,
where $P\to V$ is a $U(k+1)$-principal bundle.
We cut out the neighborhood $\C N_V$ and glue back in
$P\times _{U(k+1)}\widetilde {D_{\eps }}$, where
$\widetilde {D_{\eps }}$ is the $\eps ^{\prime }$-blow-up
of the standard 2k+2-ball of radius $\eps >\eps ^{\prime }$
(see \cite{ms1} page 250 for further details).

%\makebox(10,20)[r]{
It is easy to see that
$P\times _{U(k+1)}\widetilde {D_{\eps }}$ is
symplectomorphic to the tautological disc bundle over
$\B P_V$  the
projectivization of the normal bundle $\nu _V$.
Thus there is a natural embedding $\B P_V\to \blowMV$.
and we have the following commutative diagram in which
the horizontal arrows are symplectic embeddings.
%}
%\makebox(10,20)[l]{
$$
\xymatrix{
\cp ^k \ar[d]        & \\
\B P_V \ar[d]^{\pi } \ar[r]^f & \blowMV \ar[d]^p \\
V                    \ar[r]^i & M
}\label{D:blow-up}
$$

\centerline{\bf The blow-up diagram}

\BS

\subsection{The choice of a homology class and an
almost complex structure}

\BS

In order to define Gromov-Witten invariants we have to
chose a second
homology class  and
an almost complex structure on $\blowMV$. Moreover,
by the choice of this data we want to ensure that
some Gromov-Witten invariant is nontrivial. The idea is
to find a $J$-holomorphic stable map and to show that
this is the only element in its homology class intersecting
generic representatives of certain homology
classes.

\NI
{\bf The Choice of a homology class:}
Let $u:\cp ^1 \to \blowMV$ be the composition of the
canonical embedding of the line into $\cp ^k$ with the
maps $\cp ^k\to \B P_V \to \blowMV$ of the above
diagram. Then  we define $A:=u_*[\cp ^1]$.

\NI
{\bf The choice of an almost complex structure:}
First we define an almost complex structure
$J(\B P_V)$ on $\B P_V$ such that it is standard when
restricted to the fibers and the projection
$\pi :\B P_V\to V$ is holomorphic with respect to
compatible almost complex structure on $V$.
We decompose the tangent bundle of $\B P_V$ into the tangent bundle
along fibers $\T _{vert}\B P_V$ and the horizantal subbundle $Hor$
which is the symplectic complement of $\T _{vert}\B P_V$.  Then we
define the
almost complex structure $J(\B P_V)$ as the direct sum of
standard complex structures on $\T _{vert}\B P_V$
and the pull back of a compatible almost
complex structure on $V$ to $Hor = \pi^* \T V$.
With respect to this almost complex
structure, the projection $\pi$ becomes holomorphic.
Finally, the needed almost complex structure $J$
on $\blowMV$ is an extension of $J(\B P_V)$ compatible
with the symplectic structure.

\begin{lemma}\label{AL:nontransversal}
Let $v:\cp ^1 \to \blowMV $ be a $J$-holomorphic
curve representing the class $A=u_*[\cp ^1]$. Then the image
$v(\cp ^1)$ is contained in a fiber $\cp ^k$ of the
bundle $\pi:\B P_V\to V$.
\end{lemma}

\pf
First we
show that the image of $v$ is contained in
$\B P_V$. Start with the
observation that the intersecion pairing
$[\B P_V] \circ A = -1$.
Indeed,

\begin{eqnarray*}
[\B P_V]\circ A &=& \left <eu(\nu _{\B P_V}),A \right >\\
&=&
\left <c_1(\nu _{\B P_V}),u_*[\cp ^1]\right > \\
&=&
\left <c_1(u^*\nu _{\B P_V}), [\cp ^1]\right > = -1
\end{eqnarray*}

\NI
because the normal bundle of  $\B P_V$ restricted
to $\cp ^k\subset \B P_V$ is the tautological
line bundle.
According to
the positivity of intersection of $J$-holomorphic
submanifolds, every such curve must be
cointained in $\B P_V$.

The above agrument shows that $v:\cp ^1\to \B P_V$
is a $J(\B P_V)$-holomorphic curve.
Hence, if $v: \cp ^1 \to \B P_V$ is a $J(\B P_V)$-holomorphic map, then
$\pi \circ v: \cp ^1 \to V$ is also
holomorphic with respect to almost complex structure on $V$.
Since we know that a line contained in a fiber is
one of holomorphic curves in our homology class, $\pi \circ f$ is
null-homotopic,
especially, null-homologous. The only pseudo-holomorphic map in a symplectic
manifold, which is null-homologous, is constant.
Hence $\pi \circ v$ is a constant map, so $v$ is a map to one of fibers
of $\pi$.

\qed

\BS

\subsection{Regularity and compactness}

\BS

Let $u:\cp ^1\to (\blowMV,\om )$ be an inclusion of
the line as desribed above.
We prove the regularity of the almost complex structure
$J$ at $u$, which is stated in the following

\begin{thm}\label{AT:reg}
The linearized operator

$$D_u=D\overline{\partial }_J(u): C^{\infty }(u^*\T \blowMV)
\to \Om ^{0,1}(u^*\T \blowMV )$$

\NI
is surjective.
\end{thm}

\pf
Let $u:\cp ^1\to (\blowMV,\om )$ be as above.
There is the following splitting of the pull back
of the tangent bundle. Recal that $\dim V = 2m$

\begin{eqnarray*}
u^*(\T\blowMV )&=& u^*(\T\B P_V \oplus \nu _{\B P_V})\\
&=& u^*( \C E ^{m} \oplus \T_{vert}\B P_V \oplus \nu _{\B P_V})\\
&=&\underbrace {\C E^{m} \oplus \T\cp ^k|_{\cp ^1}}_E\oplus
\underbrace{  u^*(\nu _{\B P_V})}_N
\end{eqnarray*}

\NI
Here $\C E^n$ denotes the trivial complex vector bundle of
rank $n$, $\T_{vert}$ denotes the subbundle tangent
to the fibers of the fibration and $\nu _{\B P_V}$
is the normal bundle to $\B P_V$ in $\blowMV $.

With respect to the above splitting, the operator
$D_u:C^{\infty }(E\oplus N)\to \Om ^{0,1}(E\oplus N)$
has the following matrix form

$$
\left [
\begin{array}{cc}
D_u|_E & *    \\
0      & Du|_N
\end{array}
\right ].
$$
\BS

\NI
Now we have to show that appropriate restricted operators
are surjective.
To deal with the operator $D_u|_N$ we use the
result of Hofer, Lizan and Sikorav \cite{hls} (Theorem $1^{\prime}$ ).
It states that if $L \to \Si $ is a holomorphic line bundle
over a Riemannian surface of genus $g$, then a generalized
$\overline \p $-operator is surjective, provided that
$c_1(L) \geq 2g -1$. Thus we have to check that
$c_1(\nu _{\B P_V}) \geq -1$. But this is clearly satisfied,
because the normal bundle restricted to the fiber of
$\B P_V$ is just the tautological line bundle
over $\cp ^k$. Hence its first Chern class is -1.

To show surjectivity of $D_u|_E$ we again argue with
the direc decomposition of $E$ into $\C E^m$ and
$\T \cp^k|_{\cp ^1}$.
Since $u:\cp ^1 \to \cp ^k$ is the standard
holomorphic embedding and the structure $J$ restricted to
$\cp ^k$ is the stanard integrable one,
then the factor $\T \cp^k|_{\cp ^1}$ is preserved.
Moreover, the restriction of
the operator $D_u$ is the usual Cauchy-Riemann
derivative in this case. According to Lemma 3.5.1 in \cite{ms2} we obtain
that $D_u$ restricted to $\T \cp^k|_{\cp ^1}$ is surjective.
The surjectivity on the trivial factor $\C E^m$ again follows
from \cite{hls}.

\qed

\BS

\subsection{ The main result}

In this section we prove the nontriviality of certain
Gromov-Witten invariant of the blow-up $\blowMV$.

\begin{thm} \label{AT:main}
Let $\al,\be\in H^*(\blowMV ;\B Z )$ be cohomology classes.
Suppose that

{\em
$$\pi _*(\PD f^*(\al )) \circ \pi _*(\PD f^*(\be )) = 1 \in H_0(V).$$}
Then the Gromov-Witten invariant
$\Phi _A(\al,\be) \neq 0$.
\end{thm}

\pf
First note that

$$\PD f^*(\al ) = f_! (\PD \al ) =
[\B P_V \cap C_{\al }],$$
where
$f_!:H_*(\blowMV) \to H_{*-2}(\B P_V)$ is the homology
transfer map and $C_{\al }$ denotes a cycle representing
the homology class Poincare dual to $\al $.
Thus the latter is the homology class in $\B P_V$ obatined by
the intersection of the cycle $C_{\al }$ with
$\B P_V$.

The condition
$\pi _*(\PD f^*(\al )) \circ \pi _*(\PD f^*(\be )) = 1$
implies that for generic cycles representing the
homology classes $\PD f^*(\al )$ and $\PD f^*(\be )$
there exists an odd number of fibers, say
$\pi ^{-1}(x_1),...,\pi ^{-1}(x_{2l-1})$,
intersected by each of the cycles in exactly one
point.

Let $v:\cp ^1\to \blowMV $ be a $J$-holomorphic
curve in the class $A$ that contributes to the Gromov-Witten
invariant $\Phi _A(\al , \be )$.
Since any $J$-holomorphic curve
in the class $A$ has to be contained in some fiber
of $\pi $ (Lemma A.1), then we have that
$v(\cp ^1) \subset \pi ^{-1}(x_i)=\cp^k$ for some
$i=1,2,...,2l-1$.
Thus it follows that $\Phi _A(\al ,\be)\neq 0$, because
in each fiber $\pi ^{-1}(x_i)=\cp ^k$ there is exactly
one line passing through  two generic points and the
number of the fibers is odd.

\qed

\BS

\begin{ex}
{\em
Let $\al = a^k\in H^{2k}(\blowMV )$, where $a \in H^2(\blowMV)$ is the
class Poincar\' e dual to $[\B P_V]$. Then
$f^*(a)= c_1(\nu _{\B P_V})$ is the first Chern class of the normal
bundle of $\B P_V$.
Recall that this normal bundle when restricted to the fibre
$\cp^k$ of $\B P_V$ is the tautological line bundle
over projective space.
This implies that $f^*(a)^k = f^*(\al )$ restricts
to the positive generator in $H^{2k}(\cp ^k)$.

We claim that $\pi _*(\PD f^*(\al ) = [V]$ the fundamental
class of $V$. This follow from the following computation.
Let $\C V$ denote the positive
generator of the top cohomology of $V$.

\begin{eqnarray*}
\left <\pi _*(\PD f^*(\al ), \C V \right > &=&
\left <(\PD f^*(\al ),\pi ^* (\C V )\right >\\
&=& \left <f^*(\al )\cup \pi^*(\C V), [\B P_V] \right >\\
&=& 1.
\end{eqnarray*}

Let $\be \in H^{2n - 2}(\blowMV)$ be the class Poincare
dual to $-A\in H^2(\blowMV)$. The computation in the proof
of Lemma \ref{AL:nontransversal} shows that
$-A \cdot [\B P_V] = 1$. Note that
$\PD f^*(\be ) = f_!(-A) = 1 \in H_0(\B P_V)$ hence
$\pi _*(\PD f^*(\be )) = 1 \in H_0(V)$.
Finally we obtain that

$$\pi _*(\PD f^*(\al )) \circ \pi _*(\PD f^*(\be )) = 1 \in H_0(V)$$
and according to Theorem \ref{AT:main} we get that
$\Phi _A(\al,\be) \neq 0$.
}
\end{ex}

Now we want to give more general version of the
above example. In order to do this,
consider the following composition of maps

$$
\begin{CD}
\B P_V @>f>> \blowMV @>t>>       Th(\nu _{\B P_V}),
\end{CD}
$$

\NI
where $Th(\nu_{\B P_V})$ is the Thom space.
Let $\tau _{a } \in H^*(Th(\nu _{\B P_V}))$ be an element
corresponding to given
$a\in H^*(\B P_V)$ under the Thom isomorphism.

\begin{cor}\label{AC:main}
Let $a = \pi ^*(v) \cup c_1(\nu _{\B P_V}) ^{k-1}$ and
$b=\pi ^*(w) \cup c_1(\nu _{\B P_V}) ^{k-1}$, where
$v \cup w = \C V$ the positive generator of the top cohomology
of $V$. Then

$$\Phi _A (t^*(\tau _a),t^*(\tau _b)) \neq 0.$$
\end{cor}

\pf
First observe that
$$f^*(t^*(\tau  _{a })) = a \cup c_1(\nu _{\B P_V}),$$

\NI
because first Chern class equals the Euler class in this case.
This implies the first equality in the followin computation.

\begin{eqnarray*}
&&\pi _*(f^*(\tau _a)) \circ \pi _*(f^*(\tau _b))=\\
&&\pi _*[\PD (\pi ^*(v) \cup c_1(\nu _{\B P_V})^k)] \circ
  \pi _*[\PD (\pi ^*(w) \cup c_1(\nu _{\B P_V})^k)] = \\
&&\left <\pi ^! (\pi ^*(v) \cup c_1(\nu _{\B P_V})^k) \cup
  \pi ^! (\pi ^*(w) \cup c_1(\nu _{\B P_V})^k), [V]\right > = \\
&&\left < v \cup w , [V]\right >= 1.
\end{eqnarray*}

\NI
The second one is the definition of the transfer and
the third follows from the fact that the transfer is
a homomorphism of modules. The statement follows by
the direct application of the Theorem \ref{AT:main}.
\QED

\NI
Authors addresses: 

\NI
{\bf Jaros{\l}aw K\c edra} \\
Institute of Mathematics \\
University of Szczecin \\
ul.Wielkopolska 15 \\
70-451 Szczecin\\
Poland\\
e-mail: kedra@sus.univ.szczecin.pl

\BS
\NI
{\bf Kaoru Ono}\\
 Department of Mathematics\\
 Hokkaido University\\
Sapporo 060-0810 \\
Japan\\
 email: ono@math.sci.hokudai.ac.jp

\end{document}